\documentclass [12pt]{article}
\usepackage{amsmath}
\usepackage{amsfonts, amssymb}
\usepackage[cp1251]{inputenc} 

\newtheorem{thm}{Theorem}

\newtheorem{prop}[thm]{Proposition}

\def\R{{\mathbb R}}
\def\E{{\mathbb E}}
\def\P{{\mathbb P}}
\def\Z{{\mathbb Z}}
\def\N{{\mathbb N}}

\def\be#1\ee{\begin{equation}#1\end{equation}}
\newcommand{\bea}{\begin{eqnarray}}
\newcommand{\eea}{\end{eqnarray}}
\newcommand{\beaa}{\begin{eqnarray*}}
\newcommand{\eeaa}{\end{eqnarray*}}

\def\AA {{\mathcal{A}}}

\def\eps{{\varepsilon}}
\def\F {{\mathcal{F}}}
\def\ga{{\gamma}}

\def\M {{\mathcal{M}}\,}

\def\T{\mathbb T}

\def\x{u}

\def\tv{\widetilde V}

\begin{document}
\title{{\sc CYCLIC BEHAVIOR OF MAXIMA
\\ IN A HIERARCHICAL SUMMATION SCHEME}\footnote{Supported 
by RFBR grants 10-01-00154, 11-01-12104-ofi\_m,
and by Federal Programme 2010-1.1.-111-128-033.
}}
\author
{M.A. Lifshits}
\date{ }
\maketitle

\begin{abstract}
Cyclic behavior of maxima  in a hierarchical summation scheme.
Let i.i.d. symmetric Bernoulli random variables be associated to the edges of
a binary tree having $n$ levels. To any leaf of the tree, we associate the sum 
of variables along the path connecting the leaf with the tree root.  Let $M_n$ 
denote the maximum of all such sums. We prove that, as $n$ grows, the distributions
of $M_n$ approach some helix in the space of distributions. Each element of this 
helix is an accumulation point for the shifts of distributions of $M_n$.
\end{abstract}


\section{Branching random walks}

The study of extremal positions of branching random
walk and branching Brownian motion is by now considered as a classical problem
with first deep results obtained as early as in Hammersley's work 
\cite{Ha}. During last decade it regained a considerable popularity. New substantial 
advances were obtained but many questions remain open.

Let us shortly recall the notion of branching random walk, a very special case of which 
will be considered in this article. At initial (zero) time 
there is one particle located at zero. At time 1 the particle 
dies but gives birth to a point process (configuration) of progeny that
consists of a random number of particles (points on the real line)
whose positions are, generally speaking, mutually dependent. 
Every new born particle also lives one unit of time and dies  
giving birth to a point process of progeny independent of all other analogous 
processes. The distribution of progeny process for every particle (ancestor)
differs from the progeny process of initial particle by translation
to the position of the ancestor.
 
 Therefore, the branching random walk is a genealogical Galton--Watson tree $\T$, 
 where every element $\x\in \T$ is additionally characterized by its position on the line 
 $V(\x)$. Clearly,  $V(\x)$ is a sum (over the set of ancestors of $\x$),
 of independent random variables, each term of the sum being the displacement of a
 particle with respect to the location of its parent particle.   

We shall not consider variations of this basic model, e.g. those with random life 
time of each particle.

Every particle $\x$ belongs to a certain generation $|\x|$, i.e. to 
a level of the tree $\T$. For the initial particle, we let the generation number 
be zero.

The extremal positions in generations describing the generations' range 
are of special interest. They are defined by formulas
\[
     M_n:= \max\{V(\x), |\x|=n \} ,\qquad m_n:= \min\{V(\x), |\x|=n \}. 
\]
The limit theorems for the distributions of these variables are obtained in
\cite{ABR, Ai,Ba, BK, Br1, Br2, LS1, LS2,HuShi}. They essentially assert that 
 \be \label{limMn}
   M_n= cn+ b_n +\widetilde M_n,
\ee
where $c$ is a non-negative constant, $b_n$ is a deterministic sequence 
varying slower than a linear function  (most commonly, $b_n$ behaves 
logarithmically), and a sequence $\widetilde M_n$ converges in distribution to some limit
or is just bounded in probability. We stress that no multiplicative
norming is needed, i.e. the family of distributions of the variables 
$(M_n)_{n\ge 0}$ is shift-compact.
 
 Obviously, the linear term in the asymptotics of $M_n$ can be trivially 
eliminated by a constant shift of the progeny point process in the definition 
of the walk. 
 
 The following recent theorem due to  E. A\"{\i}d\'ekon \cite{Ai} is one of the most 
 repre\-sent\-ative and powerful results on the extremal positions.

\begin{thm} \label{t:ai} Assume that the distribution of the progeny process 
in a branching random walk is non-lattice and that the following assumptions are 
satisfied,
\be \label{ai1a}
  \E\left(\sum_{|\x|=1} 1\right)>1, 
\ee
\be \label{ai1bc}
  \E\left(\sum_{|\x|=1} e^{V(\x)} \right) =1, \qquad 
   \E\left(\sum_{|\x|=1} V(\x) e^{V(\x)} \right)=0, 
\ee
as well as the moment restrictions
\bea \label{ai2}
  && \E\left(\sum_{|\x|=1} V(\x)^2 e^{V(\x)} \right)<\infty, \qquad
  \\
  && \E\left( X(\ln_+ X)^2\right)<\infty, \qquad
  \E\left( \widetilde X(\ln_+ \widetilde X)\right)<\infty,
\eea
where $X:=\sum_{|\x|=1} e^{V(\x)}$, $\widetilde X:=\sum_{|\x|=1} V(\x)_- e^{V(\x)}$.

Then there exists an a.s. positive random variable $D$ such that for any $r\in \R$
it is true that
\[
  \lim_{n\to\infty} \P\left(M_n\le - \frac 32 \ln n +r\right) = \E\, e^{-De^{-r}}.
\]
\end{thm}

Theorem \ref{t:ai} means that in \eqref{limMn} we have $c=0$, 
$b_n= -\frac 32 \ln n$, and the distributions of $\widetilde M_n$ converge to a mixture
of shifted double exponential distributions (Gumbel laws).

Assumption \eqref{ai1a} is natural: it means that the branching process is supercritical.
This condition provides sufficient number of particles in the walk. 
Assumptions   \eqref{ai1bc} mean that a linear scaling of the walk steps 
"killing" a linear term  in \eqref{limMn} is performed.

These assumptions are not too restrictive in the following sense.
Consider a branching random walk satisfying condition
\eqref{ai1a} but not necessarily satisfying conditions
\eqref{ai1bc}. Let
\[
   \Phi(\ga):= \E\left(\sum_{|\x|=1} e^{\ga V(\x)} \right),
   \quad
   \Psi(\ga ):=\ln \Phi(\ga ), \qquad \ga >0.
\]
Making a linear shift in one generation
\[
   \tv(\x):=\ga V(\x)-\Psi(\ga ), \qquad |\x|=1,
\]
corresponds to the shift of all particles
\be \label{varch}
   \tv(\x):=\ga V(\x)-|\x|\Psi(\ga ), \qquad \x\in \T.
\ee
Let us search for $\ga >0$ such that the analogues of \eqref{ai1bc}  
for the new walk
\[
  \E\left(\sum_{|\x|=1} e^{\tv(\x)} \right) =1, \qquad 
   \E\left(\sum_{|\x|=1} \tv(\x) e^{\tv(\x)} \right)=0. 
\]
would be valid. Note that the first equality holds automatically,
since
\[
  \E\left(\sum_{|\x|=1} e^{\tv(\x)} \right) 
  =e^{-\Psi(\ga )} \E\left(\sum_{|\x|=1} e^{\ga V(\x)} \right) 
  =e^{-\Psi(\ga )}\Phi(\ga )=1.
\]
We may rewrite the second condition as
\begin{eqnarray*}
  0 &=& \E\left(\sum_{|\x|=1} (\ga V(\x)-\Psi(\ga )) e^{\ga V(\x)} \right)
  \\
  &=& \ga  \,  \E\left(\sum_{|\x|=1} V(\x) e^{\ga V(\x)} \right)
  -\Psi(\ga ) \E\left(\sum_{|\x|=1}  e^{\ga V(\x)} \right)
  \\
  &=& \ga \Phi'(\ga )-\Psi(\ga )\Phi(\ga ),
\end{eqnarray*}
which is equivalent to
\be \label{R0}
  R(\ga ):= \ga \Psi'(\ga )-\Psi(\ga ) =0.
\ee

It follows easily from H\"older inequality that the function $\Psi(\cdot)$  
is convex. Therefore, $R'(\ga )=\ga \Psi''(\ga )\ge 0$ , i.e.  $R(\cdot)$ 
is an increasing function. We have $R(0)=-\Psi(0)=-\ln\Phi(0)<0$ by \eqref{ai1a}. 
Hence, if
\[
   \Phi(\ga )<\infty, \qquad \qquad 0\le \ga <\infty,
\]
and
\be \label{rinfty}
   \lim_{\ga \to\infty} R(\ga ) 
   =  \lim_{\ga \to\infty} \left[ \ga \Psi'(\ga )-\Psi(\ga )\right] >0, 
\ee
then equation \eqref{R0} has a solution $\ga >0$, and a liner change 
\eqref{varch} reduces the study of the initial walk to the study of a walk
satisfying assumptions \eqref{ai1bc}.

In the example we focus on below, no shift can render the distributions of
$M_n$ convergent to a non-degenerate limit distribution. Instead,
they approach some helix of distributions, or, if the shifts are allowed, 
they circulate along some closed curve in the space of distributions. 
There are {\it two} reasons preventing application of Theorem \ref{t:ai} 
in that case: first, the distributions of the progeny process is a lattice 
one; second, the reduction condition \eqref{rinfty} fails.

\section{Hierarchical summation scheme}

In the following we consider the simplest model of a branching random
walk: every particle produces {\it two} particles whose translations are
independent Bernoulli random variables taking value $1$ with probability
$p$ and $-1$ with probability $1-p$. Therefore, the genealogic tree
$\T$ is just a simplest binary tree and the particles' locations are described
by the sums of independent Bernoulli random variables along the branches of
this tree. It is amazing that such a simple model demonstrates an interesting
limit behavior.

We may redescribe the object under consideration as follows.

Consider $n$-level binary tree and associate to its edges i.i.d. Bernoulli 
random variables $(B_i)$. The tree has $2^n$ leafs. To each leaf we associate 
the sum of random variables picked up along the path connecting the leaf with 
the tree root. Let $M_n$ be the maximum of the sums along all leafs. We shall 
investigate the asymptotic behavior of the distribution of $M_n$, as $n$
goes to infinity. Clearly, we have, $M_0=0$, $M_n\in[-n,n]$, and  $
M_n=n\ (\textrm{mod}\ 2)$. Moreover, there is a recurrency equation 
 \be 
    M_{n+1}=  \max\left\{ M_n^{(1)}+B^{(1)}; M_n^{(2)}+B^{(2)}\right\},
 \ee    
 where $M_n^{(j)}$ and $B^{(j)}$ are independent copies of $M_n$, resp.
 of the Bernoulli variable. 

It is worthwhile to notice that hierarchical summation schemes appear not only 
in connection to the branching walks. They emerge, for example, in physical models 
such as Derrida generalized random energy studied by Bovier and Kurkova \cite{BK}. 
In their setting, the summands situated on the same level of the tree have the
same distribution but this distribution is allowed to vary reasonably from one 
level to another. 

\subsection{Symmetric case}

In this subsection we consider the most interesting symmetric case
 \[
     \P(B_i=1)=\P(B_i=-1)= \frac 12\ .
 \]
 
 We start with the study of the behavior of $\E M_n$. Subsequent 
 delicate considerations are entirely based on the following modest fact.
 
\begin{prop} \label{prop1} 
Let 
\[
  K_n:=\{u: |u|=n, V(u)=M_n\}
\]
be the number of vertices of level $n$ where the maximum $M_n$ is attained.
Then $K_n\to\infty$ in probability and
\[
  \lim_{n\to \infty} \E(M_{n+1}-M_n) =1.
\]
\end{prop}

{\bf Proof.}\ Notice that $K_n$ is bounded from below by a critical
Galton--Watson process $Z_n$ with the progeny number $N$ defined by the law
\[ 
  \P(N=k)=\begin{cases} \tfrac 14, &k=0, \\
                         \tfrac 12,&k=1, \\
                         \tfrac 14,&k=2,
           \end{cases}                
\]
and restarting from 1 at extinction time. To observe $Z_n$ on the tree, it is 
enough to keep track of the paths along which we have only  $B_i=+1$; at each 
level, where we have extinction (the values $-1$ occupy all continuations of  
the paths we observe), we keep a single path and consider only its continuations -- 
according to the previous rule. Remark that all chosen paths  provide maximal 
values of sums on each level, hence $Z_n\le K_n$. 

Look at $Z_n$ from the point of view of Markov chain theory. All states are
recurrent and null, since the expectation of extinction time for our 
Galton--Watson process is infinite. Hence, for any fixed $\ell\in\N$
\[
   \lim_{n\to\infty} \P(Z_n=\ell)=0,
\]
e.g. see Theorem 3 in \cite[Section XIII.3]{Fel}. Hence, for any $m\in \N$
\[
  \lim_{n\to \infty} \P(K_n\le m) \le \lim_{n\to \infty} \P(Z_n\le m)
  = \sum_{\ell=1}^m \lim_{n\to \infty} \P(Z_n =\ell) =0, 
\]
as claimed.

Passing to the expectations, let us notice that
$M_{n+1}-M_n\in\{-1,+1\}$; moreover,
\begin{eqnarray*}
 \P(M_{n+1}-M_n=+1|\, \AA_n) &=& 1-2^{-K_n},
 \\
 \P(M_{n+1}-M_n=-1|\, \AA_n) &=& 2^{-K_n},
\end{eqnarray*}
where  $\AA_n$ stands for the sigma-field generated by the variables situated on
first $n$ levels of the tree. It follows that
\[
  \E(M_{n+1}-M_n)= 1- 2\, \E\,  2^{-K_n},
\]
and the second claim of the proposition follows from the first one.
$\Box$
\medskip 

Proposition \ref{prop1} shows that $\E M_n\sim n$, as $n$ grows to infinity. 
Hence, it suggests that $M_n$ is relatively close to its
upper border $n$. Therefore, it is more convenient to consider
the variables $M_n'=\frac{n-M_n}2$. Then $M_n'$ is a non-negative  
integer random variable and satisfies the relations $M_0'=0$, 
$M_n'\in[0,n]$ and the equation
 \be 
     M_{n+1}'= \min
     \left\{ M_n'^{(1)}+\tilde B^{(1)}; M_n'^{(2)}+\tilde B^{(2)}\right\}
 \ee         
where $M_n'^{(j)}$ and $\tilde B^{(j)}$ are independent copies of $M_n'$, 
resp. of  a variable $\tilde B$ having the distribution
\[
     \P(\tilde B=1)=\P(\tilde B=0)= \frac 12\ .
 \]
It is more convenient to express the recurrency equation in terms of the 
tails of random variables. Let $F_n(x):=\P(M_n'\ge x)$. Then
\[ 
    F_0(x)=\begin{cases} 1,& x\le 0,
                 \\ 0,& x>0,
            \end{cases}
\]
and
\be \label{itera}
   F_{n+1}(x)=\left[\frac{F_n(x)+F_{n}(x-1)} 2 \right]^2.
\ee
This equation has {\it many} invariant solutions. Indeed, 
an invariant solution should satisfy equations
\be \label{invar}
  4 F(x)=\left[F(x)+F(x-1) \right]^2.
\ee
Hence, $F(x-1)= G(F(x))$ and $F(x)=g(F(x-1))$, where 
$G(y):= 2\sqrt{y}-y$ and $g(y):=2-y-2\sqrt{1-y}$ are mutually inverse
functions.
It follows that all values of $F$ can be expressed via $F(0)$ 
as iterations of functions $g$ and $G$.
The family of invariant distribution may be written in a parametric 
form $\{\F^a, {0 < a < 1}\}$, where
\[
    \F^a(n)= \begin{cases}
    g^n(a), & n>0,\\
    a,      & n=0,\\
    G^{|n|}(a),& n<0.
    \end{cases}      
\]
and $g^n, G^n$ denote the $n$-th iteration of $g$, resp. $G$. It
is clear that the family of invariant distributions form a
continuous one-parametric curve (it is natural to call it a "helix") 
in the space of distributions $\M(\R^1)$; moreover, using the appropriate
shifts we can transform this curve into a closed cycle, i.e. 
$\F^{g(a)}(\cdot-1)=\F^a(\cdot)$ for any $0<a<1$.

We consider now the limit behavior of $F_n(x)$ as $n\to \infty$.  
Let us first handle the case of fixed $x$. The following is true.

\begin{prop} \label{prop2} 
   For each $x\in \Z$ it is true that $\lim_{n\to\infty} F_n(x)=1$.
\end{prop}

{\bf Proof}. 
Using induction in $x$, we derive from \eqref{itera} that the sequence 
$F_n(x)$ is non-decreasing in $n$ for each fixed $x$. 
Hence, the limit $F(x):=\lim_n F_n(x)$ exists and satisfies
equation \eqref{invar}. Notice that $F(x)=1$ for $x\le 0$,
and it follows from \eqref{invar} that $F(x-1)=1$ yields $F(x)=1$. 
Therefore, $F(x)=1$ for each $x\in \Z$.
$\Box$  

It is worthwhile to notice that the fist non-trivial case, $x=1$, 
corresponds to the behavior of $\P(M_n<n)$, i.e to the extinction
probability of a critical branching Galton--Watson process. It is well 
known from the classical works of R. Fisher and A.N. Kolmogorov that
for such process extinction takes place almost surely. 

Proposition \ref{prop2} implies that the variables $M_n'$ converge 
to infinity in probability.

We pass now to the main result called a cyclic limit theorem. 
We will show that for large $n$ the distribution $F_n$ admits an 
approximation by an appropriate invariant distribution.

\begin{thm} \label{limsim}
For each $n$ let a median $k_n$ be defined by the relation 
\[
   k_n=\inf\{x\in \Z : F_n(x)\le 1/2\}.
\]
Let $a_n=G^{k_n}(F_n(k_n))$. Then
\[
    \lim_{n\to\infty} \sup_{x\in \Z} | F_n(x)-\F^{a_n}(x)| = 0.
\]
\end{thm}

Since  $F_n$ goes along the limit helix $(\F^a)_{0<a<1}$ with decreasing 
speed, it is natural to conjecture that all points of the helix are the limit 
points of $F_n$ after appropriate shift normalizations. In particular no
shifts can render $F_n$ convergent to a non-degenerate distribution. 
The exact assertion is as follows. 

\begin{thm} \label{limpoint}
For any $a\in (0,1)$ there exist $z\in \Z$ and a sequence 
of integers  tending to infinity $n_k$ such that  
\[
    \lim_{k\to\infty} \max_{x\in \Z} |\F^a(x)- F_{n_k}(x+k-z)| =0.  
\]
\end{thm}
\medskip

{\bf Proof of Theorem \ref{limsim}.} \ 

First, notice that the second claim of Proposition  \ref{prop1} may be 
also written in the form
\be \label{Delta}
   \lim_{n\to \infty} \Delta_n :=  \lim_{n\to \infty} \E(M'_{n+1}-M'_n) 
   = \lim_{n\to \infty} \sum_{k=1}^\infty [F_{n+1}(k)-F_n(k)] =0.
\ee

Another necessary ingredient for the proof is the identity
\be \label{identgF}
  F_n(x)=g(F_n(x-1)-\delta)-\delta,
\ee
where $\delta=\delta(n,x):=F_{n+1}(x)-F_{n}(x)\in [0,\Delta_n]$.

Indeed, we may write \eqref{itera} as 
\[
  F_n(x)+\delta= \left(\frac{F_n(x)+\delta + F_n(x-1)-\delta}{2} \right)^2.
\]
Taking into account that the function $g$ satisfies the identity
$g(y)=\left(\tfrac{g(y)+y}{2} \right)^2$,
we arrive at \eqref{identgF}. Note for subsequent applications that \eqref{identgF}
yields useful inequalities
\be \label{FngG}
   F_n(x)\le g(F_n(x-1)), \quad  G(F_n(x))\le F_n(x-1).
\ee 

Now we pass to the proof of the theorem.
By using \eqref{identgF} and the monotonicity of $g(\cdot)$, for  
each integer $d\ge 0$ we obtain
\begin{eqnarray*}
  F_n(k_n+d) &=& g(F_n(k_n+d-1)-\delta) -\delta 
  \\
  &\le& g(F_n(k_n+d-1))\le ... \le g^d(F_n(k_n)).
\end{eqnarray*}
We also have
\[
  \F^{a_n}(k_n+d) = g^d( \F^{a_n}(k_n))   =g^d(F_n(k_n)).
\]
For any $\eps>0$, take a positive integer $D$ such that $g^D(\frac 12)\le \eps$. 
Then for any $d\ge D$, by using monotonicity of $g(\cdot)$ 
and inequality $g(y)\le y$, we infer
\[
   \max\{  F_n(k_n+d); \F^{a_n}(k_n+d)\}
   \le g^d(F_n(k_n)) \le g^d(\tfrac 12) \le g^D(\tfrac 12)\le \eps. 
\]
Hence,
\be \label {dD1}
   \max_{d\ge D} | F_n(k_n+d)- \F^{a_n}(k_n+d)| \le \eps. 
\ee
Now we show by induction that for every $d=0,1,\dots, D$
it is true that
\be   \label {dD2}
   \lim_{n\to \infty} | F_n(k_n+d)- \F^{a_n}(k_n+d)| =0.
\ee
We have chosen parameters $a_n$ so that 
\[
   \F^{a_n}(k_n) = g^{k_n}(\F^{a_n}(0)) = g^{k_n}(a_n)
   = g^{k_n}G^{k_n} (F_n(k_n))= F_n(k_n). 
\]   
Therefore, for $d=0$ the left hand side of \eqref{dD2} vanishes, thus
providing the induction base.  Assume that for $d-1$ assertion \eqref{dD2} 
is proved, then by \eqref{itera} for $d$ we have
\[
  | F_n(k_n+d)- \F^{a_n}(k_n+d)| 
  = | g(F_n(k_n+d-1)-\delta) - g(\F^{a_n}(k_n+d-1))| +\delta, 
\]
where $\delta:= F_{n+1}(k_n+d)- F_{n}(k_n+d)\in [0,\Delta_n]$.
It follows that
\begin{eqnarray*}
  &&| F_n(k_n+d)- \F^{a_n}(k_n+d)| 
  \\
  &\le&  \left[ | F_n(k_n+d-1) - \F^{a_n}(k_n+d-1)|+\Delta_n\right] 
       \max_{0\le y\le \frac 12}|g'(y)| +\Delta_n,
\end{eqnarray*}
Since $\Delta_n \to 0$ by \eqref{Delta}, and since the function $g'$ 
is bounded on $[0,\frac 12]$, we obtain that
\begin{eqnarray*}
   && \limsup_{n\to \infty} | F_n(k_n+d)- \F^{a_n}(k_n+d)|
   \\
   &\le&
  \lim_{n\to \infty} | F_n(k_n+d-1)- \F^{a_n}(k_n+d-1)|\cdot 
  \max_{0\le y\le \frac 12}|g'(y)| 
  =0.
\end{eqnarray*}
Therefore, \eqref{dD2} is proved. By combining \eqref{dD1} with \eqref{dD2},
we obtain
\[
    \lim_{n\to \infty} \max_{d\ge 0} | F_n(k_n+d)- \F^{a_n}(k_n+d)| =0.
\]
Negative $d$'s are handled in the same way by using function $G$ instead 
of $g$. $\Box$
\medskip

{\bf Proof of Theorem \ref{limpoint}.} \ 
Without loss of generality we may assume that 
$\tfrac{1}{2}\not\in \left\{\F^a(x),x\in\Z\right\}$. Then there exists
$z\in \Z$ such that
\[  
   \F^a(z-1)>\frac 12 >\F^a(z).
\]
Fix an $\eps>0$ and choose $\delta\in (0,\min\{a,1-a\})$ so small that
$b\in (a-\delta,a+\delta)$ implies
\[
   \max_{x\in \Z} |\F^b(x)-\F^a(x)|<\eps.
\]
We may also require the inequalities 
\be \label{alsomed}
   \F^{a-\delta}(z-1)> \frac 12 > \F^{a+\delta}(z).
\ee  
to hold. Take a positive integer $n_0$ such that for all $n\ge n_0$ it is
true that
$\Delta_n < \F^{a+\delta}(z)- \F^{a-\delta}(z)$.
Let now $k$ be so large that $F_{n_0}(k)< \F^{a-\delta}(z)$.
Consider the sequence $f_n:= F_n(k), n\ge n_0,$ for fixed $k$.
By Proposition \ref{prop2} it grows to one. Since  
$f_{n_0}< \F^{a-\delta}(z)$ and for all $n\ge n_0$ it is true that 
\[ 
  f_{n+1}-f_n=  F_{n+1}(k)-F_n(k)\le \Delta_n
  \le \F^{a+\delta}(z)- \F^{a-\delta}(z),
\]
there exists $n:=n_k$ satisfying
\[
   F_n(k)=f_n\in (\F^{a-\delta}(z), \F^{a+\delta}(z)).
\]
Notice that $k$ is the median for $F_n$, since by \eqref{alsomed}, \eqref{FngG}
\[
  F_n(k)\le \F^{a+\delta}(z)<\frac 12\,
  ;
  F_n(k-1)\ge G(F_n(k))\ge G(\F^{a-\delta}(z)) =\F^{a-\delta}(z-1)> \frac 12.
\]
Therefore, the approximating distribution $\F^{a_n}$ from Theorem \ref{limsim}
satisfies the equalities
\[
   \F^{a_n}(k)=F_n(k)=\F^b(z)
\]
for some $b\in (a-\delta,a+\delta)$. Finally, we use the following
fact: if $\F^a(u)=\F^b(v)$ for some
$a,b\in (0,1)$ and some $x,y\in \Z$, then for all $x\in \Z$ we have
\[
    \F^a(x+u-v)= \F^b(x).
\]
In our case $\F^{a_n}(x+k-z)= \F^b(x)$ holds. Therefore,
\begin{eqnarray*}
&& \max_{x\in \Z} |\F^a(x)- F_n(x+k-z)|
\\
&\le& 
\max_{x\in \Z} |\F^a(x)- \F^b(x)| + \max_{x\in \Z} |\F^b(x)- F_n(x+k-z)| 
\\
&\le&
\eps + \max_{x\in \Z} |\F^{a_n}(x+k-z)- F_n(x+k-z)|.
\end{eqnarray*}
Since $\eps$ was chosen arbitrarily and the second term tends to zero by 
Theorem \ref{limsim}, we obtain the assertion of Theorem \ref{limpoint}.\
$\Box$
\medskip

One of the reasons for non-existence of the unique limit distribution
is the discrete type of Bernoulli distribution, as is clearly seen
from Theorem \ref{t:ai}. Another, less obvious and may be a deeper,
reason is the failure of \eqref{rinfty}. 
Indeed, in the hierarchical summation scheme for $p$-Bernoulli variables 
we have
\begin{eqnarray*}
  \ga\ \Psi'(\ga)-\Psi(\ga) 
  &=& \ga\ \frac{pe^\ga-(1-p)e^{-\ga}}{pe^\ga+(1-p)e^{-\ga}} 
      -\ln 2-\ln\left(pe^\ga+(1-p)e^{-\ga}\right)
   \\   
   &=& -\ln (2p)- \frac{2(1-p)\ga}{pe^{2\ga}}(1+o(1))
\end{eqnarray*}
and
\[
   \lim_{\ga\to\infty} \left[  \ga\Psi'(\ga)-\Psi(\ga) \right] =-\ln (2p). 
\]
Therefore, condition \eqref{rinfty} is satisfied iff $p<\tfrac 12$.

The tree structure of the hierarchical summation scheme is not related to
the helix-type behavior of the maxima distributions: one can obtain a
similar result for conventional summation (see Section \ref{s:notree} 
below).
\medskip

{\bf Remark.} One can also derive Theorem \ref{limsim} from Theorem 1 in
Bramson's work \cite{Br1}. The additional advantages of his result are the more
general branching rule and approximation in the sense of almost sure convergence. 
However, Theorem \ref{limsim} provides more transparent geometric picture of 
the phenomenon.

\subsection{A limit theorem for the case $p>1/2$} \label{s:Blimit}

We maintain the notation of the previous subsection but assume now that 
 \[
    \P(B_i=1)=1-\P(B_i=-1)= p >1/2 \ .
 \]
 Let $q:=1-p$. The recurrency equation now takes the form
 \be 
    M_{n+1}'= \min\left\{ M_n'^{(1)}+\tilde B^{(1)}; M_n'^{(2)}+\tilde B^{(2)}\right\}
 \ee         
where $M_n'^{(j)}$ and $\tilde B^{(j)}$ stand for independent copies of $M_n'$ and of
a variable $\tilde B$ that satisfies 
\[
    \P(\tilde B=1)=1-\P(\tilde B=0)= q.
\]
In terms of the distribution tails $F_n(x):=\P(M_n'\ge x)$, we obtain an equation 
analogous to (\ref{itera}), namely,
\be \label{itera2}
   F_{n+1}(x)=\left[ F_n(x)p+F_{n}(x-1)q  \right]^2.
\ee
There is a big difference with respect to the previous case: now there exists 
a {\it unique} invariant non-degenerate solution satisfying the equation
\be \label{invar2}
   F(x)=\left[F(x)p+F(x-1)q \right]^2
\ee
and the initial condition $F(x)=1, x\le 0$. Namely, 
\[
F(x)= (2p^2)^{-1} \left[ 
1-2F(x-1)pq - \sqrt{1-4F(x-1)pq}\right],
\qquad x>0.
\]
Therefore, it is not surprising that a limit theorem holds in this case.

\begin{thm} Uniformly over $x\in \Z$, the monotone convergence 
$F_n(x)\nearrow F(x)$ holds.
\end{thm}

{\bf Proof.}\  First, by induction in $x$ we derive from \eqref{itera2} that 
the sequence $F_n(x)$ is non-decreasing in $n$ for each fixed $x$. 
Therefore, the limit $F(x):=\lim_n F_n(x)$ exists and satisfies equation 
\eqref{invar2}. It remains to prove that it is non-degenerate, i.e. it is
different from identical unit. For this purpose, it is enough to notice that
\[
   F_n(1)=\P(M_n'\ge 1)=\P(n-M_n\ge 2)=\P(M_n\le n-2)
\]
coincides with extinction probability of the {\it supercritical} 
Galton-Watson process with the progeny number $N$ defined by the law
\[ 
  \P(N=k)=\begin{cases} q^2, &k=0, \\
                         2pq,&k=1, \\
                         p^2,&k=2.
           \end{cases}                
\]
Therefore, $1-F(1)$ is the survival  probability of the process, which is 
strictly positive, as $p>\frac{1}{2}$.
\ $\Box$

\subsection{Some results for the case $p< 1/2$} \label{s:Bdrift}

In what concerns limit theorems, not more is known for this case than for the
hierarchical summation scheme with general independent random variables 
having finite exponential moments.
For $p=P(B=1)<1/2$ the equation (\ref{invar2}) has no non-trivial solutions, 
therefore, the behavior of maxima is completely different than in the previous cases
-- a drift with constant speed appears. Once we eliminate this linear drift, the 
distributions of $M_n$ form a dense set with exponentially decreasing tails.

Recall that a family of random variables $(X_n)$ is called
{\it shift-compact}, if there exists a real sequence $(a_n)$ such that the 
distributions of random variables $X_n-a_n$ form a tight family on the real line, 
i.e.
\[ 
    \lim_{K\to\infty} \sup_n \P\{|X_n-a_n|>K \}=0.
\]
\medskip

\begin{prop} \label{prop3} 
Let $p<1/2$. Then the sequence of random variables $M_n$ is shift-compact, 
while
\[ 
   \E M_n \sim \rho \, n, \qquad  n \to\infty,
\] 
where the shift coefficient $\rho$ is defined from equation
\be \label{Bdrift}
     2p^\rho q^{1-\rho}= \rho^\rho(1-\rho)^{1-\rho}.
\ee 
\end{prop}

{\bf Proof.}\ The result follows, e.g., from Theorem 1.1 in  \cite{BZ}. 
It is worthwhile to notice that the equation for the drift \eqref{Bdrift} 
is essentially the special case of equation \eqref{R0} providing reduction 
to the critical case. 
$\Box$
\bigskip

\section{Cyclic theorems for maxima of independent sums}
\label{s:notree}

Let  $(\xi_i)_{i\in \N}$ be {\it integer}\ i.i.d. random variables.   
Consider the sum $S_n:=\sum_{i=1}^n \xi_i$,  and let
$S_n^{(j)}$, $1\le j\le 2^n$, be independent copies of $S_n$. 
We are interested in the behavior of $M_n:=\max_{j\le 2^n} S_n^{(j)}$. 

We will assume that our random variables satisfy
\be \label{moments}
    \E|\xi_1|<\infty\quad \textrm{and}\quad \E\exp\{\ga \xi_1\}<\infty,\ \forall \ga>0.
\ee

Let $\omega$ be the upper bound of the distribution,
\[ 
   \omega:= \sup \{m\in \N: \P(\xi_1=m)>0 \}. 
\]
Assume that one of the two following assumptions is satisfied: either

$(i)$ \qquad $\omega=\infty$,

\noindent or

$(ii)$ \qquad  $\omega<\infty$\ and\ $\P\left(\xi_1=\omega\right)<1/2$. 
\medskip

Since the cumulant
\[
   L(\ga):=\ln \E \exp\{\ga \xi_1\}
\]   
is a convex function of $\ga$, the function $L(\ga)-\ga L'(\ga)$ 
is non-increasing. It is continuous and vanishes at $\gamma=0$.
Moreover, if \eqref{moments} holds, and any 
of assumptions $(i)$  or $(ii)$ is satisfied, it is easy to show that 
\[ 
    \lim_{\ga\to +\infty} [ L(\ga)-\ga L'(\ga)] < \ln (1/2).
\]
Therefore, a solution of equation
\be \label{gammastar}
    L(\ga)-\ga L'(\ga) = \ln (1/2)
\ee 
exists on $(0,+\infty)$.  Let denote it $\ga_*$ and let
$\rho_*:= L'(\ga_*)$. Notice also that under either $(i)$ or 
$(ii)$ the distribution of $\xi_i$ is non-degenerated (not concentrated at a
single point), therefore the solution of \eqref{gammastar} is unique.

\begin{thm} \label{cycle_notree} Let \eqref{moments}  and either 
$(i)$  or $(ii)$ holds. Let  
$\rho_*,\gamma_*$ be defined by equation $(\ref{gammastar})$. Then  
\be \label{bound_gen}
   \P\left\{ M_n < \rho_* n - \frac{\ln n}{2 \gamma_*} +z \right\}
   =\exp\left\{ - \, \frac{\exp\{-\gamma_* z\} (1+o(1))}
    {\sqrt{2\pi}\sigma(\ga_*)(1-e^{-\ga_*})}  \right\},
\ee
where $\sigma(\cdot)^2=L''(\cdot)$, uniformly
over\footnote{In other words,
 we consider $z$ such that the expression in the left hand side is an 
 integer number.}  
\[
    z\in I\bigcap  \left[\Z - \rho_* n  + \frac{\ln n}{2 \gamma_*} \right]
\]   
 for any bounded interval $I$.  
\end{thm}

We can rewrite formula \eqref{bound_gen} as
\be \label{bound_geni}
   \P\left\{ M_n < m \right\}
   =\exp\left\{ - \exp\{-\gamma_* (m-a_n)\} (1+o(1)) \right\},
   \qquad m\in \Z,
\ee
where
\[
  a_n:= \rho_* n   
        -\frac{\ln[\sqrt{2\pi n}\sigma(\ga_*)(1-e^{-\ga_*})]}{\ga_*}\ .
\]

For each $a\in\R$ let $\F^a$ denote the distribution on integers
given by
\[
   \F^a((m,+\infty))
   = \exp\left\{ - \exp\{-\gamma_* (m-a)\}\right\}, \qquad m\in\Z.
\]
Then $(\F^a)_{a\in\R}$ is a curve in the space of distributions. It is natural
to perceive it as a helix, in view of 1-periodicity up to a shift: 
$\F^{a+1}\{m+1\}=\F^{a}\{m\}$. Relation \eqref{bound_geni} shows that the 
distribution of r.v. $M_n$ admits the uniform approximation by the helix
element $\F^{a_n}$, while after appropriate centering it admits an approximation
by the element  $\F^{[a_n]}$ of the helix turn $(\F^a)_{0\le a < 1}$.
Moreover any distribution $(\F^a)_{0\le a < 1}$ is a limit of some subsequence
of centered distributions of $M_n$.

The proof of Theorem \ref{cycle_notree}, which is supposed to be published separately,
is based on a large deviation theorem due to V.V. Petrov \cite[Complement 2 in \S4 
Chapter VIII]{Petr}.

Let us consider Bernoulli case as an example.
Let $\xi_i=B_i$ be independent random variables having non-symmetric 
Bernoulli distribution, i.e.
\[
   \P(B_i=1)=1-\P(B_i=-1)= p < 1/2 \ .
 \]
Let the drift coefficient $\rho_*$ be again defined by equation (\ref{Bdrift}). 
 We also need two auxiliary constants 
$\kappa:=\frac{p(1-\rho_*)}{q\rho_*}\in (0,1)$ and $\beta:=2\pi\rho_*(1-\rho_*)$.
 Then the result of Theorem \ref{cycle_notree} takes the following form.

 \begin{thm} \label{cycle_notree_b} 
  We have
 \be \label{bound_ber}
     \P\left\{ M_n< \rho_* n -\frac{\ln(\beta n)}{2|\ln\kappa|} +z \right\}
     =\exp\left\{ - \frac{\kappa^{z}}{1-\kappa}\, (1+o(1))\right\},
 \ee
 uniformly over
 \[
    z\in I\bigcap  \left[\Z - \rho_* n +\frac{\ln(\beta n)}{2|\ln\kappa|}\right]
 \]   
 for any bounded interval $I$.
 \end{thm}

{\bf Remark.} For $p \ge \tfrac 12$ neither of conditions $(i), (ii)$ holds.
Equation \eqref{gammastar} has no solutions, thus Theorem \ref{cycle_notree} 
does not apply. 
\bigskip

The author is deeply indebted to Irina Kurkova and to Zhan Shi for interesting discussions,
for providing important references, and, most of all, for motivation to write this article.

\bibliographystyle{amsplain}

\begin{thebibliography}{10}

{\baselineskip=12pt

\bibitem{ABR} L. Addario-Berry, B. Reed, 
{\it Minima in branching random walks}. -- Ann. Probab. {\bf 37} (2009), 
1044--1079.

\bibitem{Ai} E. A\"{\i}d\'ekon,   
{\it Convergence in law of the minimum of a branching random walk}. -- Ann. 
Probab. (to appear). Preprint {\tt arxiv:\,1101.1810}\ (2011).

\bibitem{Ba} M. Bachmann,   
{\it Limit theorems for the minimal position in a branching random walk with 
independent logconcave displacements}. -- Adv. Appl. Probab. {\bf 32} (2010), 
159--176.

\bibitem{BK} A. Bovier, I. Kurkova, 
{\it Derrida's generalized random energy models. I.: models with finitely many 
hierarchies}. -- Ann. Inst. H. Poincar\'e, {\bf 40} (2004), 439--480.

\bibitem{Br1} M. Bramson, 
{\it Minimal displacement of branching random walk}. -- Z. Wahrsch. Theor., 
{\bf 45} (1978), 89--108.

\bibitem{Br2} M. Bramson,   
{\it Convergence of solutions of the Kolmogorov equation to travelling waves}. 
-- Mem. Amer. Math. Soc.  {\bf 44} (1983), No 285.

\bibitem{BZ} M. Bramson, O. Zeitouni, 
{\it Tightness for a family of recursive equations}. Ann. Probab., 
{\bf 37} (2009), 615--653.

\bibitem{Fel} W. Feller,   
{\it An Introduction to Probability Theory and Its Applications}. Vol.I,
2-nd ed., Wiley, N.Y., 1957. 

\bibitem{Ha} J.M.Hammersley,  
{\it Postulates for subadditive processes}.-- Ann. Probab., 
{\bf 2} (1974), 652--680.

\bibitem{HuShi} Y. Hu, Z. Shi,   
{\it Minimal position and critical martingale convergence in branching random 
walks, and directed polimers on disordered trees}. -- Ann. Probab. {\bf 37} 
(2009), 742--789.

\bibitem{LS1} S.P. Lalley, T. Selke, 
{\it A conditional limit theorem for the frontier of the branching Brownian 
motion}. --  Ann. Probab. {\bf 15} (1983), 1052--1061.

\bibitem{LS2} S.P. Lalley, T. Selke, 
{\it Limit theorems for the frontier of a one-dimensional branching motion}. -- 
Ann. Probab. {\bf 20} (1992), 1310--1340.

\bibitem{Petr} V.V. Petrov,   
{\it Sums of Independent Random Variables}. Nauka, Moscow, 1972. 

 }
\end{thebibliography}

St.Petersburg State University

email: {\tt  lifts@mail.rcom.ru}
\pagebreak

\end{document}